\documentclass{amsart}
\usepackage{color}
\usepackage{amsmath}
\usepackage{amssymb}
\usepackage{amsfonts}

\newtheorem{Proposition}{Proposition}

  \newtheorem{Lemma}[Proposition]{Lemma}
  
  \newtheorem{Theorem}[Proposition]{Theorem}

 \newtheorem{Note}[Proposition]{Note}

\newcommand {\z}{{\noindent}}
\def\blackslug{\hbox{\hskip 1pt \vrule width 4pt height 8pt depth 1.5pt
\hskip 1pt}}
\def\qed{\quad\blackslug\lower 8.5pt\null\par}

 \def\RR{\mathbb{R}}
 \def\NN{\mathbb{N}}

\begin{document}

\title[Stokes constants]{Rigorous bounds of Stokes constants for some
  nonlinear ODEs at rank one irregular singularities}

\author{Ovidiu Costin, Rodica D. Costin and Matthew Kohut}

\address{Department of Mathematics, Hill Center, Busch Campus,
        Rutgers, The State University of New Jersey,  
       110 Frelinghuysen Rd, Piscataway, NJ
        08854-8019, costin\symbol{64}math.rutgers.edu}
\maketitle

\today

\begin{abstract}
  
  A rigorous way to obtain sharp bounds for Stokes constants is
  introduced and illustrated on a concrete problem arising in
  applications.

\z {\bf Keywords}: Stokes constants, exponential asymptotics, resurgence.
\end{abstract}

\section{Introduction}
Stokes constants (multipliers) relate to the change in asymptotic
behavior of a solution of a differential equation as the direction
toward an irregular singularity is changed (cf. \S\ref{Sec2}).  If the
constants are nonzero, then the change in behavior of the solution is
nontrivial and this fact plays a very important role in a number of
problems.

Many interesting results are known for linear ordinary differential
equations, see \cite{Balser, Immink, Loday, Lutz} and references
therein; papers \cite{Daal-Jap, OldeDaalhuis} use hyperasymptotic
methods to express Stokes constants as convergent series.

Stokes multipliers have been evaluated in closed form for a
wide class of {\em integrable systems} sometimes using difficult and
subtle arguments, see e.g. \cite{Deift, Flashka, Its, Miwa, Miwa2} and
references therein. Integrability however is non-generic, and is
thought to play a crucial role in any {\em explicit} evaluation of
Stokes constants.

A theory of wide applicability of extended Borel summation, Borel
plane singularities and their relation witgh Stokes phenomena and was
introduced by \'Ecalle, \cite{Ecalle-book}. For generic nonlinear
systems complete asymptotic expansions (transseries) of solutions and
their Borel summability, singularities in Borel plane and formulas
linking them to nonlinear Stokes transitions are rigorously obtained
in \cite{Duke}.  The paper \cite{OC-MDK-PLond} finds the link between
the structure of singularities of solutions in the Borel plane and
summation to the least term, as well as with the behavior of the
coefficients of the asymptotic series in generic nonlinear systems.

In applications it often only matters whether the Stokes constants are
nonzero, while their exact value is not relevant. The results in
\cite{Duke} and \cite{OC-MDK-PLond} provide a rather straightforward
way of obtaining rigorous and accurate {\em estimates} of the Stokes
multipliers for a large class of linear or nonlinear differential
systems; in principle any prescribed precision can be obtained, as
well as the information that a constant does not vanish if such is the
case.

One goal of the present paper is to present such a method and at the
same time complete an argument in the proof by Tanveer and Xie
\cite{Tanveer-Xie} for the nonexistence of steady fingers with width
less than $\frac{1}{2}$ when small nonzero surface tension is taken
into account. Their argument relies on a conjectured nonzero value of a Stokes
constant of the differential equation
\begin{equation}\label{eqv}
2v'' - t +\frac{1}{v^2} = 0
\end{equation}
Eq.
(\ref{eqv}) appears as an ``inner-equation'' arising in the context of
steady Hele-Shaw cell fingers (see also \cite{CDHPP}). It convenient to
illustrate our general technique through this particular equation.  We
show that two Stokes constants for (\ref{eqv}) are given by
\begin{equation}\label{S1}
S_1=ib\frac{\pi^{3/2}2^{13/14}}{\Gamma(1/7)\Gamma(3/7)}
\end{equation}
\begin{equation}\label{S2}
S_2=i{\rm{e}}^{\frac{i\pi}{14} }b\frac{\pi^{3/2}2^{13/14}}{\Gamma(1/7)\Gamma(3/7)}
\end{equation}
with 
\begin{equation}\label{estimb0}
 1\leq b \leq 1+\frac{12}{37}
\end{equation}
In particular it follows from (\ref{estimb0}) that indeed $S_1$ and
$S_2$ are nonzero\footnote{Numerical calculation gives
  $b=1.1722\cdots$}.

It will become apparent that the method introduced here applies to
generic systems of equations whose irregular singularity has rank one
(this is the most frequent type of irregular singularities in
applications). It relies on a detailed relation established in
\cite{OC-MDK-PLond} between the Stokes constants and the behavior for
large index of the coefficients of the asymptotic series solutions,
followed by inductive proof of bounds on the solution of the
recurrence relation defining them.

\section{Stokes constants and exponential asymptotics}\label{Sec2}
Consider a system of differential system of the form
\begin{equation}\label{eqy}
{\mathbf{y}}'
=\left(-\mathbf{\Lambda}-
  \frac{1}{x} \mathbf{B}\right){\mathbf{y}}
+\mathbf{f_0}(x)+\mathbf{g}(x,{\mathbf{y}})
\end{equation}
where $\mathbf{\Lambda}={\rm{diag}}\boldsymbol{\lambda}$,
$\boldsymbol{\lambda}=(\lambda_1,...,\lambda_n)$,
$\mathbf{B}={\rm{diag}}\boldsymbol{\beta}$,
$\boldsymbol{\beta}=(\beta_1,...,\beta_n)$, $\mathbf{f_0}=O(x^{-2})$,
$\mathbf{g}=O(x^{-2};|{\mathbf{y}}|^{2};|x^{-2}{\mathbf{y}}|)$ (as
$|x|\rightarrow\infty$), with $\mathbf{f_0}$ analytic at $\infty$ on a
half-line $d$ and $\mathbf{g}$ analytic at $(\infty,0)$ under {\em
  nonresonance assumptions} \cite{Duke} (a slightly weaker condition
than the linear independence of the eigenvalues $\lambda_j$ over the
rationals). The system (\ref{eqy}) has then a {\em rank one irregular
  singularity} at infinity.
\subsection{Power series solutions and exponentially small terms}

The general type of formal solutions of differential systems in the
presence of irregular singular points was studied in detail by Fabry
\cite{Fabry} and Cope \cite{Cope} (see also \cite{Duke},
\cite{OC-MDK-PLond}, \cite{Invent}). For (\ref{eqy}) (assumed
nonresonant) the general formal solution for large $x$ has the form
\begin{equation}\label{genformsol}
\tilde{\mathbf y}=\tilde{\mathbf
  y}_0+\sum_{\mathbf{k}\in\NN^n\setminus \mathbf{0}}\mathbf{C}^\mathbf{k}{\rm{e}}^{-\mathbf{k}\cdot\boldsymbol{\lambda}x}x^{\mathbf{k}\cdot\boldsymbol{\beta}}\tilde{\mathbf{s}}_\mathbf{k}
\end{equation}
where $\tilde{\mathbf y}_0$ and $\tilde{\mathbf{s}}_\mathbf{k}$ are
power series (generically divergent) and $\mathbf{C}=(C_1,...,C_n)$
are free parameters.

The general formal solutions (\ref{genformsol}) can in fact be
calculated algorithmically, in a way that that will be briefly
sketched.

The power series solution $\tilde{\mathbf y}_0$ is unique. To
determine terms beyond all orders of this series, a formal calculation
is carried out to find formal solutions which are small perturbations
to $\tilde{\mathbf y}_0$; for instance one substitutes
$\tilde{\mathbf{y}}=\tilde{\mathbf{y}}_0+\boldsymbol{\delta}$ in
(\ref{eqy}), where $\boldsymbol{\delta}\ll x^{-q}$ for all $q\in\RR_+$
(as $x\rightarrow \infty$ on some direction). Since
$\tilde{\mathbf{y}}_0$ satisfies (\ref{eqy}) we formally get
$$\boldsymbol{\delta}'\sim\left(-\mathbf{\Lambda}-
  \frac{1}{x} \mathbf{B}\right)\boldsymbol{\delta}$$
hence in a first approximation
\begin{equation}\label{delta0}
\boldsymbol{\delta}\sim
{\rm{e}}^{-\mathbf{\Lambda}x}x^{-\mathbf{B}}\mathbf{C}
\end{equation}
where $\mathbf{C}$ is a vector of free parameters.

Order by order perturbation expansion around (\ref{delta0}) first
produces power series multiplying the exponential (\ref{delta0}) and
then smaller and smaller exponentials, eventually leading to
(\ref{genformsol}).
\subsection{Transseries solutions}
From the point of view of correspondence of formal solutions to actual
solutions it was recognized that in general only {\em asymptotic}
expressions of the form (\ref{genformsol}) ({\em{transseries}}, as
introduced by \' Ecalle \cite{Ecalle-book}) can be lifted to actual
functions.

Transseries and their correspondence with functions constitute the
subject of exponential asymptotics, a field which developed
substantially in the eighties with the work of Berry
(hyperasymptotics), \' Ecalle (the theory of analyzable functions) and
Kruskal (tower representations and nice functions).

Let $d_\theta={\rm{e}}^{i\theta}\RR_+$ be a direction in the complex
$x$ plane. A {\em{transseries solution along $d_\theta$}} is, in our
context, a formal solution (\ref{genformsol}) whose terms are well
ordered with respect to the relation: $\gg$ as $x\rightarrow\infty$,
$x\in d_\theta$. In particular, a formal solution (\ref{genformsol})
is a transseries along $d_\theta$ if constants $C_j=0$ for all $j$
such that ${\rm{e}}^{-\lambda_jx}$ is not going to zero for $x\in
d_\theta$.
\subsection{Exact solutions associated to transseries solutions}
For any direction $d_\theta$ there exists a one-to-one correspondence
between transseries solutions along $d_\theta$ and actual solutions
that go to zero on this direction. The correspondence is natural,
constructive and compatible with all operations with functions
(respectively, transseries) \cite{Ecalle-book}.

Actual solutions $\mathbf y$ corresponding to a transseries
(\ref{genformsol}) on a direction $d_\theta$ have the same classical
asymptotic expansion $\tilde{\mathbf y}_0$ as $x\rightarrow\infty$,
$x\in d_\theta$. The constants $C_j$ multiply small exponentials in
the transasymptotic expansion of the solutions $\mathbf y$, and are
beyond all orders of the power series $\tilde{\mathbf y}_0$. They
therefore cannot be defined using the classical Poincar\'e definition
of asymptoticity.

For nonresonant systems with a rank 1 irregular singular point this
correspondence was established in \cite{IMRN} and \cite{Duke}. It was
shown that the series $\tilde{\mathbf y}_0$ is Borel summable (in a
generalized sense); its inverse Laplace transform ${\mathbf Y}_0(p)$
is analytic at $p=0$ and is (generically) singular at an array of
points (which determine the Stokes directions). If we denote
\begin{equation}\label{denotey0}
\tilde{\mathbf y}_0(x)\, =\, \sum_{n=2}^{\infty}\, \frac{1}{x^n}\, {\mathbf
  y}_{0;n}
\end{equation}
then we have
\begin{equation}\label{denoY0}{\mathbf Y}_0(p)\, =\, \sum_{n=2}^{\infty}\, \,
\frac{p^{n-1}}{\Gamma(n)}\, {\mathbf y}_{0;n}
\end{equation}

The other power series $\tilde{\mathbf{s}}_\mathbf{k}$
($|\mathbf{k}|>0$) in a transseries solution (\ref{genformsol}) are
also Borel summable in the generalized sense. It is then shown in
\cite{Duke} that the series of functions obtained from a transseries
solution by replacing the power series with their Borel sums is
convergent for $x$ large enough on appropriate directions. The
solution thus obtained is asymptotic to the power series
$\tilde{\mathbf y}_0$ for large $x$ on the direction on which
summation was performed. The correspondence obtained between
transseries solutions and actual solutions is one-to-one and is
compatible with all operations.

\subsection{Stokes phenomenon}

Consider a solution $\mathbf y$ that goes to zero for
$x\rightarrow\infty$ on a direction $d_\theta$. To associate a 
transseries $\tilde{\mathbf y}$ to $\mathbf y$ along $d_\theta$ means
in fact to specify the parameters $C_j$, hence $C_j$ may be different
for different values of $\theta$.

It turns out, however, that the $C_j$, as functions of $\theta$, are piecewise 
constant; the directions $d_\theta$ at which one of the $C_j$ has
a jump discontinuity are called {\em{Stokes directions}} (see
\cite{Duke} for more details).

The way Stokes multipliers $S$ are related to the classical asymptotic
behavior of the solutions is given in the following Proposition
of \cite{Duke}:

\begin{Proposition}\label{T01}  Consider eq. (\ref{eqy}) under the
  assumptions given. Assume $\lambda_1$ is an eigenvalue of least
  modulus. Without loss of generality we can assume that $\lambda_1=1$
  and $\Re\beta_j<0$ for all $j$ (these inequalities can be
  arranged in the normalization process, \cite{Duke}). Let $\gamma^{\pm}$ be two
  paths in the right half plane, near the positive/ negative imaginary
  axis such that $|x^{-\beta_1+1}e^{-x\lambda_1}|\rightarrow 1$ as
  $x\rightarrow\infty$ along $\gamma^{\pm}$. Consider the solutions
  $\mathbf{y}$ of (\ref{eqy}) which are small in any proper subsector
  of the right half plane. 
  
  Then, along $\gamma^{\pm}$ we have, for some $C$,
\begin{gather}\label{classicS}
  \mathbf{y}= (C\pm\frac{1}{2}S_1)\mathbf{e}_1
  x^{-\beta_1+1}e^{-x\lambda_1}+o(\mathbf{e}_1
  x^{-\beta_1+1}e^{-x\lambda_1})
\end{gather}
for large $x$ along $\gamma^{\pm}$ (where $\mathbf{e}_j$ is the $j^{th}$
unit vector $\mathbf{e}_j=(0,...,0,{\underbrace{1}}_{j},0,...,0)$).
\end{Proposition}

Stokes constants relate to the Maclaurin series of $\mathbf{Y}_0$ in
the following way (\cite{Duke} and \cite{OC-MDK-PLond}\footnote{There
  is a typo in formula (2.1) of \cite{OC-MDK-PLond}: $\beta'$ should
  read $\beta$ (as obtained in its proof).}):
\begin{Theorem}[\cite{OC-MDK-PLond}]\label{T1} Under the same assumptions as
  in Proposition \ref{T01} we have
\begin{equation}\label{mainth}
\mathbf{Y}^{(r)}_0(0)\, =\, \sum_{j;|\lambda_j|=1}\,
\frac{\Gamma(r-\beta_j+1)}{2\pi i {\rm{e}}^{i(r+1-\beta_j)\phi_j}}\, \left(\,
S_j\mathbf{e}_j+\mathbf{h}_j(r)\, \right)
\end{equation}
where $\mathbf{Y}_0(p)$ is the generalized inverse Laplace transform
of $\tilde{\mathbf{y}}_0$ (see (\ref{denotey0}), (\ref{denoY0})),
$\phi_j=\arg \lambda_j$ (ordered increasingly, starting with
$\lambda_1=1$, $\phi_1=0$), $\mathbf{h}_j(r)\sim
r^{-1}\mathbf{h}_{j;0}$ for large $r$.
\end{Theorem}
\section{Main Results}
{\bf Normalization}. To apply Theorem \ref{T1} to solutions of
(\ref{eqv}), the equation has to be shown to be amenable to the normal
form (\ref{eqy}). The substitution
\begin{equation}\label{normsub}
v(t)\, =\, t^{-1/2}\, \left( \, 1+\, u(x)\, \right)\ \ \ \
{\mbox{where}}\, \ x=\frac{4}{7}\, t^{7/4}
\end{equation}
transforms (\ref{eqv}) to 
\begin{equation}\label{equ}
u''=u+\frac{1}{7}\, \frac{1}{x}\, u'-\frac{12}{49}\, \frac{1}{x^2}\,
u-\frac{3u^2+2u^3}{2(1+u)^2}\, -\frac{12}{49}\, \frac{1}{x^2}
\end{equation}
Substitution (\ref{normsub}) is natural and a general procedure for
finding normalizing substitutions was described in \cite{Invent}.

Equation (\ref{equ}) can be written as a system
\begin{equation}\label{sysequ}
{\mathbf{u}}'
=\left(\begin{array}{cc} 0 & 1 \\ 1 & \frac{1}{7x}  \end{array}\right){\mathbf{u}}
+\left(\begin{array}{c} 0\\ -\frac{12}{49}\end{array}\right)\frac{1}{x^2}+
\left(\begin{array}{l} 0\\ h(x,\mathbf{u})\end{array}\right) 
\end{equation}
where
$${\mathbf{u}}=\left(\begin{array}{l} u_1\\u_2\end{array}\right)\ \ \ 
,\ \ \ h(x,\mathbf{u})=-\frac{12}{49}\, \frac{1}{x^2}\,
u_1-\frac{3u_1^2+2u_1^3}{2(1+u_1)^2}$$
The dominant linear part of
(\ref{sysequ}) is diagonalized by substituting
\begin{equation}\label{subuy}
{\mathbf{u}}(x)=S(x){\mathbf{y}}(x)\ \,\ \ {\mbox{with}}\ \ S(x)=\left(\begin{array}{cc} 1 & 1 \\ -1+\frac{1}{14x} &
    1+\frac{1}{14x} \end{array}\right)
\end{equation}
which gives the normal form (\ref{eqy}) 
with $n=2$ and
\begin{equation}\label{LamdaB}\mathbf{\Lambda}=\left(\begin{array}{cc} 1 & 0 \\ 0 & -1
  \end{array}\right)\ \ \ ,\ \ \ \ \mathbf{B}=\left(\begin{array}{cc} -\frac{1}{14} & 0 \\ 0 & -\frac{1}{14}
  \end{array}\right) 
\end{equation}
\begin{equation}\label{f0}
\mathbf{f_0}(x)=\frac{1}{x^2}\mathbf{f_0}\ \ , {\mbox{with}}\ \mathbf{f_0}=\frac{6}{49} \left(\begin{array}{c} 1\\
    -1\end{array}\right)
\end{equation}
\begin{equation}\label{hereisg}
\mathbf{g}(x,{\mathbf{y}})=\frac{15}{392}\frac{1}{x^2}\left(\begin{array}{cc} -1 & -1 \\ 1 & 1
  \end{array}\right) \mathbf{y}+ \frac{1}{2}  h\left( x,S\mathbf{y}\right) \left(\begin{array}{c} -1\\
   1 \end{array}\right)
\end{equation}
Equation (\ref{eqy}) with (\ref{LamdaB}), (\ref{f0}), (\ref{hereisg})
has a rank 1 irregular singularity at $x=\infty$ and it is
written in normal form.
The eigenvalues of the matrix $\mathbf{\Lambda}$ are
$\boldsymbol{\lambda}=(1,-1)$ and $\beta_1=\beta_2=-1/14$.

The dominant power series in the transseries solution
$\tilde{\mathbf{y}}_0=(\tilde{y}_{0;1},\tilde{y}_{0;2})$ is
$$\tilde{y}_{0;1+\sigma }(x)={\frac
  {6}{49}}\,{x}^{-2}+(-1)^\sigma{\frac {87}{343}}\,{x}^{-3}+{ \frac
  {2028}{2401}}\,{x}^{-4}+(-1)^\sigma{\frac {57798}{16807}}\,{x}^{-5}+...$$($\sigma\in\{0,1\}$).

{\em{Transseries solutions}}.
For eq. (\ref{eqy}), (\ref{LamdaB}--\ref{hereisg}) the transseries solutions along $d_\theta$ with
$\theta\in(-\pi/2,\pi/2)$ must have $C_2=0$ (while $C_1$ is
arbitrary); transseries along directions with
$\theta\in(\pi/2,3\pi/2)$ must have $C_1=0$ (and $C_2$ is arbitrary).

\z {\em{Stokes phenomena.}} The problem of physical interest depends on
the solutions ${\mathbf y}$ which are classically asymptotic to the
power series $\tilde{\mathbf y}_0$ for $|x|\rightarrow\infty$ in the
sector $\arg x\in [0,5\pi/8]$ \cite{Tanveer-Xie}. The fact that the
asymptoticity is required on a large enough sector implies that this
solution is unique and its transseries on any direction of argument in
$(0,5\pi/8]$ are expansions (\ref{genformsol}) with $C_1=C_2=0$.

Indeed, a transseries (\ref{genformsol}) for $\arg x \in[0,\pi/2]$
must have $C_2=0$. Since $i\RR_+$ is not a Stokes direction however, $C_2$
remains zero in all directions with $\arg x \in (\pi/2,5\pi/8]$; but
in these directions a transseries also has $C_1=0$ and since $i\RR_+$
is not a Stokes direction, $C_1=0$ also for
$\arg x \in (0,\pi/2]$.

Since $\RR_+$ is a Stokes direction, $C_1$ may become nonzero here,
and its value is the Stokes constant $S_1$.

Relation (\ref{mainth}) for equation (\ref{eqy}),
(\ref{LamdaB}-\ref{hereisg}) is
\begin{equation}\label{FormS}
\mathbf{Y}^{(r)}_0(0)\, =\, 
\frac{\Gamma\left( r+\frac{15}{14}\right)}{2\pi i}\, \left(\,
  S_1\mathbf{e}_1\, +\, 
{S_2} {{\rm{e}}^{-i\left( r+\frac{15}{14}\right)\pi}}\,
\mathbf{e}_2+\mathbf{h}(r)\, \right)
\end{equation}
with $\mathbf{h}(r)\sim r^{-1}\mathbf{h}_{0}$ for large $r$.

From (\ref{subuy}) we have $u(x)=y_1(x)+y_2(x)$. Hence, using
(\ref{FormS}), we get:
\begin{equation}\label{FormU}
U^{(r)}_0(0)\, =\, Y^{(r)}_{0;1}(0)+Y^{(r)}_{0;2}(0)\, =
\frac{\Gamma\left( r+\frac{15}{14}\right)}{2\pi i}\, \left(\,
  S_1+
{S_2} {{\rm{e}}^{-i\left( r+\frac{15}{14}\right)\pi}}\,
+O(r^{-1})\, \right)
\end{equation}
Though obvious, it is worth pointing out that the inverse Laplace
transform $U(p)$ of $u(x)$ has the convergent series expansion at
$p=0$
\begin{equation}\label{serUp}
U_0(p)=\sum_{n\geq 1}\frac{u_{2n}}{\Gamma({2n})}p^{{2n}-1}
\end{equation}
hence $U^{(r)}_0(0)=0$ if $r$ is even and $U^{(r)}_0(0)=u_{r+1}>0$ if
$r$ is odd.
Then from (\ref{FormU}) for $r$ even it follows that 
\begin{equation}\label{S2S1}
S_2=-S_1{{\rm{e}}^{\frac{15}{14} i\pi}}
\end{equation}
Using (\ref{FormU}) for $r=2n-1$, (\ref{S2S1}) and (\ref{asyun}) we
get the formulas (\ref{S1}) and (\ref{S2}) where $b$ satisfies
(\ref{estimb}).
\section{Estimating the constant $b$}\label{Scalcb}
The constant $b$ in (\ref{asyun}) below can be estimated within any
prescribed accuracy by the procedure described in this section. For
the problem at hand however, the inequalities (\ref{estimb0}) obtained
here are more than enough.
\begin{Proposition}\label{mainprop}
Any solution of equation (\ref{equ}) that goes to zero along $\RR_+$
has the power series expansion
\begin{equation}\label{seru}
u(x)\sim\sum_{n\geq 1}x^{-2n}u_{2n}\ \ \ \ (x\rightarrow +\infty)
\end{equation}
where $u_{2n}>0$ and
\begin{equation}\label{asyun}
u_{2n}=b\frac{\sqrt{\pi}}{\Gamma(1/7)\Gamma(3/7)}n^{-13/14}\Gamma(2n+1)\left(1+O(n^{-1})\right)
\end{equation}
with $b$ satisfying
\begin{equation}\label{estimb}
1 \leq b\leq 1.324
\end{equation}
\end{Proposition}

$ $

Solutions of (\ref{eqv}) satisfying $v(t)\sim t^{-1/2}$
as $t\rightarrow+\infty$ have the asymptotic expansion
\begin{equation}\label{eq2}
v(t)=\sum_{k=0}^\infty \frac{c_k}{t^{\frac{7}{2}k+\frac{1}{2}}}
\ ,\  c_0=1
\end{equation}
where the $c_n$ satisfy the recurrence relation:
$$c_n=\frac{(7n-6)(7n-4)}{4}c_{n-1}\hskip 6cm$$
\begin{equation}\label{4}
+\frac{1}{2}\sum_{k=1}^{n-1}\frac{(7n-7k-6)(7n-7k-4)}{2}c_{n-k-1}\sum_{i=0}^kc_ic_{k-i}-\frac{1}{2}\sum_{k=1}^{n-1}c_kc_{n-k}
\end{equation}
Let $d_n$ be the sequence satisfying the recurrence 
\begin{equation}\label{5}
d_n=\frac{(7n-6)(7n-4)}{4}d_{n-1}\ \ \ \ (n\geq 1),\ \ \ d_0=1
\end{equation}
Clearly, $d_n>0$ and in fact
\begin{equation}\label{6}
d_n=\left(\frac{49}{4}\right)^n\frac{\Gamma(n+\frac{1}{7})\Gamma(n+\frac{3}{7})}{\Gamma(\frac{1}{7})\Gamma(\frac{3}{7})}
\end{equation}
Denote
\begin{equation}\label{notbn}
b_n=\frac{c_n}{d_n}
\end{equation}
From (\ref{4}), (\ref{5}) the recurrence for $b_n$ is
\begin{equation}\label{recbn}
b_n=b_{n-1}+Q_n\ \ \ ,\ \ b_0=1
\end{equation}
where $Q_n=Q_n^+-Q_n^-$ {with}
\begin{equation}\label{defQnp}
Q_n^+=\frac{1}{2d_n}\sum_{k=1}^{n-1}\frac{(7n-7k-6)(7n-7k-4)}{2}b_{n-k-1}c_{n-k-1}T_k
\end{equation}
\begin{equation}\label{defQnm}Q_n^-=\frac{1}{2d_n}T'_{n}\ \  ,\ \
  \ \ T_k=2d_kb_k+T'_{k}\end{equation}
\begin{equation}\label{defTpk}
  T'_{k}=\sum_{i=1}^{k-1}b_ib_{k-i}d_id_{k-i}
\end{equation}
\begin{Proposition}\label{L3} 
The sequence $b_n$ converges and its limit $b$ satisfies
the estimate  (\ref{estimb}).
Therefore
\begin{equation}\label{asycn}
c_n\sim  \Gamma(2n+1)\left(\frac{49}{16}\right)^n
n^{-\frac{13}{14}}\frac{b\sqrt{\pi}}{\Gamma(\frac{1}{7})\Gamma(\frac{3}{7})}
\ \ \ \ \ (n\rightarrow +\infty)
\end{equation}
\end{Proposition}
The proof of Proposition \ref{L3} relies on the following two Lemmas:
\begin{Lemma}\label{L1}
Let $n \ge 5$. Assume there exist $A_1,A_2>0$ such that 
$$A_1 \le b_k \le A_2\ \ \ {{for\ all\ }}\ k \ \ with \ \ 0\leq k\leq n-1$$
Then 
\begin{equation}\label{estimQn}
|Q_n| \le \frac{B}{n^2}
\end{equation}
where
\begin{equation}\label{defB}
B\, =\, 0.6\,  A_2^2+ 0.0144\, A_2^3
\end{equation}
\end{Lemma}
The proof of Lemma \ref{L1} is given in \S\ref{ProofL1}.
\begin{Lemma}\label{L2} 
For $A_1=1$ and $A_2=1.324$ we have
\begin{equation}\label{b0..b7}
A_1\leq b_k<A_2\ \ for\ k=0,1,2,...,7
\end{equation}
and 
\begin{equation}\label{bgeq8}
A_1+\frac{B}{k}\le b_k\le A_2-\frac{B}{k}\ \  where\ B=1.0787
,\ for\  all\ k\ge 8
\end{equation}
\end{Lemma}
The proof of Lemma \ref{L2} is given in \S\ref{ProofL2}.

\z {\em{Proof of Proposition \ref{L3}.}} This is an immediate
consequence of Lemmas \ref{L1} and \ref{L2}. Indeed, by Lemma \ref{L2}
we have $1\leq b_k<1.324$ for all $k\geq 0$. Then by Lemma \ref{L1} we
have $|Q_n|<Bn^{-2}$ for all $n\geq 5$, and by (\ref{recbn}) the
sequence $b_n$ is Cauchy. The estimate (\ref{estimb}) follows from
Lemma \ref{L2} and the asymptotic behavior (\ref{asycn}) of $c_n$
follows from (\ref{notbn}) and (\ref{6}). Relation (\ref{asycn})
follows from (\ref{6}) and the Stirling formula.

It only remains to prove Lemmas \ref{L1} and \ref{L2}.
\subsection{{Proof of Lemma \ref{L1}.}}\label{ProofL1}  
Note the following estimate for any $N\geq N_0\geq 3$:
$$\frac{1}{d_N}\sum_{i=1}^{N-1}d_id_{N-i}=2\frac{d_1d_{N-1}}{d_N}+\frac{1}{d_N}\sum_{i=2}^{N-2}d_id_{N-i}$$
$$\leq
2\frac{d_1d_{N-1}}{d_N}+(N-3)\frac{d_2d_{N-2}}{d_N}
\leq \frac{1}{N^2}E(N)\leq \frac{1}{N^2}E(N_0)\ \ \ {\mbox{for}}\
N\geq N_0\geq 3
$$
where
\begin{equation}\label{defE}
E(N)=\frac{6}{49}\frac{N^2}{(N-\frac{6}{7})(N-\frac{4}{7})}+\frac{240}{49^2}\frac{N^2}{(N-\frac{4}{7})(N-\frac{11}{7})(N-\frac{13}{7})} 
\end{equation}
Therefore
\begin{equation}\label{estimdN}
D_N\equiv\frac{1}{d_N}\sum_{i=1}^{N-1}d_id_{N-i}\leq \frac{1}{N^2} E(N_0)\ \ \ {\mbox{for}}\
N\geq N_0\geq 3 
\end{equation}
It is easy to check that the estimate also holds for $N=N_0=2$. To
estimate $T'_k$ (see (\ref{defTpk})) note that since it was assumed
that $b_k\geq A_1>0$ for all $k\leq n-1$, then for $N_0\leq {k}\leq
n-1$, we have, using (\ref{estimdN}),
\begin{equation}\label{aav}
0<T'_k \le A_2^2D_kd_k\le A_2^2E(N_0)\frac{1}{k^2}d_k \ \
\ {\mbox{for\ }}{N_0\leq k}\leq n-1
\end{equation}
Let $N_0=5$; we have $E(5)<0.24$.  For $k=2,3,4$ to estimate
(\ref{aav}) further note that
$D_kd_k=k^2D_k\frac{d_k}{k^2}<0.22\frac{d_k}{k^2}$. Then from
(\ref{aav}) we get
\begin{equation}\label{aa}
0<T'_k < 0.24\, A_2^2 \frac{1}{k^2}d_k \leq 0.06 \,  A_2^2\, d_k \ \
\ {\mbox{for\ }}{2\leq k}\leq n-1
\end{equation}
It follows that (see (\ref{defQnm}))
\begin{equation}\label{dst}
0<T_k \le \alpha d_k, \ \ {\mbox{for}}\ 1\leq k\leq
n-1\ \ ,\ \ \text{where $\alpha=2A_2+ 0.06\, A_2^2 $}
\end{equation}
For $1\le k \le {n-1}$ we have, using (\ref{dst}), (\ref{6})
\begin{equation}\label{evterm1}
0<\frac{(7n-7k-6)(7n-7k-4)}{2}b_{n-k-1}d_{n-k-1}T_k
\end{equation}
$$\leq  A_2\frac{(7n-7k-6)(7n-7k-4)}{2}d_{n-k-1}T_k$$
\begin{equation}\label{evterm2}
=2A_2d_{n-k}T_k \le 2\alpha
A_2 d_kd_{n-k}
\end{equation} 
and using (\ref{estimdN}) (see (\ref{defQnp}))
\begin{equation}\label{estQp}
0<Q_n^+\le \frac{\alpha A_2}{d_n}\sum_{k=1}^{n-1}d_kd_{n-k} =\alpha
A_2D_n\leq \frac{\alpha A_2E(N_0)}{n^2}< \frac{0.24\, \alpha A_2}{n^2}
\end{equation} 
for $n\geq 5$. The term $Q_n^-$ (see (\ref{defQnm})) is estimated
similarly, using (\ref{aa}):
\begin{equation}\label{estQm}
0< Q_n^-\leq \frac{A_2^2E(N_0)}{2n^2}<\frac{0.12 \, A_2^2}{n^2}\
{\mbox{for}}\ n\geq 5
\end{equation} 
Then (\ref{estQp}), (\ref{estQm}) implies (\ref{estimQn}) which proves
Lemma \ref{L1}.
\subsection{{Proof of Lemma \ref{L2}.}}\label{ProofL2}   A direct calculation yields:
\begin{multline*}
  b_{{1}}=1\ ,\ b_{{2}}={\frac {169}{160}}= 1.061\cdots\ ,\ b_{{3}}={
  \frac {743}{680}}= 1.092\cdots \ ,\ b_{{4}}={\frac
  {426573}{382976}}= 1.113\cdots\\
b_{{5}}={\frac
{71300607}{63289600}}= 1.126\cdots\ ,\ b_{{6}}={\frac {1406520669011}{
1239463526400}}= 1.134\cdots\\
b_{{7}}={\frac { 135335882622883}{118668949344000}}= 1.140\cdots \ 
,\ b_{{8}}={\frac {6575066918153233021 }{5744440195153920000}}=
1.144\cdots 
\end{multline*}
Then (\ref{b0..b7}) holds, and also (\ref{bgeq8}) is true for $n=8$.
Estimate (\ref{bgeq8}) is shown by induction.

Let $n\geq 9$. Assuming (\ref{bgeq8}) for all $k$ with $8\leq k\leq n-1$, then in
particular $A_1\leq b_k\leq A_2$ for all $k$ with $0\leq k\leq n-1$
hence (\ref{estimQn}) holds. Using (\ref{recbn}) in (\ref{bgeq8}) for
$k=n-1$ we get
$$A_1+\frac{B}{n-1}+Q_n\leq b_n\leq A_2-\frac{B}{n-1}+Q_n$$
which in view of (\ref{estimQn}) implies (\ref{bgeq8}) for $k=n$.
Lemma \ref{L2} is proved. 
\begin{Note}
  Substantially sharper estimates of $b$ can be obtained using more
  terms in the expansion of $\mathbf{Y}^{(r)}_0(0) $ that can be
  easily obtained from \cite{OC-MDK-PLond}.
\end{Note}
\z {\bf{Acknowledgments}} The authors would like to thank Prof. S
Tanveer for suggesting the problem. Work of RDC and MK was partially
supported by the Rutgers REU program.  MK would like to thank C
Carpenter and Profs. C Woodward and I Blank. Work of OC was partially
supported by NSF grants 0103807 and 0100495.

\end{document}